\documentclass[12pt]{article}
\usepackage[namelimits,intlimits]{amsmath}
\usepackage[english]{babel}
\usepackage{amsfonts}
\usepackage{euscript}
\usepackage{amsthm}
\usepackage{multicol,amsthm,syntonly}
\usepackage{euscript}
\usepackage{color}
\usepackage{graphicx}
\usepackage{eufrak}
\begin{document}
\begin{center}

\vspace{30mm} {\Large \bf Linear quasigroups. I}\\
\vspace{10mm} {\bf Tabarov Abdullo\footnote[1]{The research to
this article was sponsored by Special Projects Office, Special
and Extension Programs of the Central European University Corporation. Grant Application Number:321/20090.}}\\
{\it Tajik National University, Department of Mechanics and
Mathematics,\\
734017, Dushanbe, Rudaki ave.17, Tajikistan}\\
e-mail: tabarov2010@gmail.com
\end{center}

\vspace{20mm}
\pagestyle{plain}
\setcounter{page}{1}
\sloppy
\begin{center}
{\bf Abstract}
\end{center}

The article is devoted to linear quasigroups and some of their
generalizations. In the first part main definitions and notions of
the theory of  quasigroups are given. In the second part some
elementary properties of linear quasigroups and their
generalizations are presented. Finally in the third part
endomorphisms and endotopies of linear quasigroups and their
generalizations are investigated.
\\



\vspace{20mm}
\begin{center}
{\bf I. Main definitions and  notions}
\end{center}

{\bf Definition 1.1.} A groupoid  $(Q, \cdot)$ is called \emph {a
quasigroup} if the equations
$$
 a \cdot
x=b, \quad  y \cdot a=b
$$
are uniquely solvable for any $a, b \in Q$ [1].

A quasigroup can also be defined in another way - this is an
algebra $(Q, \cdot, /,  \setminus )$  with three binary operations
$( \cdot )$, $( / )$ and $( \setminus )$, satisfying the following
identities:
$$
(xy)/y=x, \,\,(x/y)y=x,\,\, y(y \backslash x)=x,\,\, y \backslash
(yx)=x.
$$

The quasigroup $(Q,\cdot)$ is called {\it isotopic} to the
quasigroup $(Q,\circ)$, if there exist three permutations $ \alpha
,\beta ,\gamma $ of the set $Q$ such that $\gamma (x \circ y) =
\alpha x \cdot \beta y$ for any  $x,y \in Q$.

{\bf Definition 1.2.} {\it The element  $f(x)$ $(e(x))$ of a
quasigroup $(Q, \cdot)$ is called left (right) local identity
element of an element   $x \in Q,$ if  $f(x)\cdot x=x$ $(x\cdot
e(x)=x)$.}

If  $f(x)\cdot y=y$ $(y\cdot e(x)=y)$ for all $y\in Q,$ then
$f(x)$ $(e(x))$  is called {\it left (right) identity element} of
a quasigroup $(Q, \cdot).$

If $f(x)=e(x)$ for all $x \in Q$ i.e. when all left  and right
local identities  of the quasigroup coincide, then  $(Q, \cdot)$
has identity element which is denoted as  $e: \,\, e(x)=f(x)=e.$ A
quasigroup with identity is called \emph{a loop}.

{\bf Theorem 1.1 (A.A. Albert [1])}. {\it Every quasigroup is
isotopic to some loop.}

An isotopy of the form $T= (\alpha, \beta,\varepsilon)$ where
$\varepsilon$ is the identity permutation is called {\it principal
isotopy}.

{\bf Theorem 1.2 (A.A. Albert [1])}. {\it If a loop $(Q, \circ)$
is isotopic to a group $(Q,+)$, then $(Q, \circ)$ is a group and
$(Q, \circ) \cong (Q,+)$, that is the groups $(Q, \circ)$ and
$(Q,+)$ are isomorphic.}

{\bf Theorem 1.3 (V.D. Belousov [1])}. {\it If the loop
$(Q,\circ)$ is principally isotopic to the quasigroup $(Q,
\cdot)$, then the isotopy must have the form: $T=(R^{-1}_{a},
L^{-1}_{b}, \varepsilon)$,} where $ L_{a} x = ax, \,\, R_{a} x =
xa, \,\, R^{-1}_{a}x= a/x,\,\, L^{-1}_{b}x=b \backslash x, $ for
all $a,b,x \in Q$.

An isotopy of the form $T=(R_{a}^{-1}, L_{b}^{-1}, \varepsilon)$
is also called  LP-isotopy of the quasigroup $(Q,\cdot)$. So any
LP-isotopy of a quasigroup is a loop.

In the class of quasigroups that are isotopic to groups the
interesting object are the so  called {\it linear quasigroups}
that were first introduced by  V.D. Belousov in [2] in connection
with researching    balanced identities in quasigroups.

{\bf Definition 1.3.} {\it  A quasigroup $(Q,\cdot)$ is called
{\it linear} over a group $(Q, + )$, if it has the form:
\begin{equation}
\label{eq1} xy = \varphi x + c + \psi y,
\end{equation}
where $\varphi ,\psi \in Aut(Q, + )$,  $c$ is a fixed element from
$Q$.}

Later G.B. Belyavskaya and A.KH. Tabarov in [3], by analogy with
linear quasigroups, defined alinear quasigroups, as well as
introduced the classes of left and right linear quasigroups, left
and right alinear quasigroups and mixed type of linearity.

{\bf Definition 1.4.} {\it A quasigroup $(Q,\cdot)$ is called {\it
alinear} over a group $(Q, + )$, if it has the form:
\begin{equation}
 xy=\bar{\varphi}x+c+\bar{\psi}y,
\end{equation}
where $\bar{\varphi}, \bar{\psi}$ are antiautomorphisms  of $(Q, +
)$, $c$ is a fixed element from $Q$.}

{\bf Definition 1.5.} {\it A quasigroup $(Q,\cdot)$ is called {\it
left (right) linear} over a  group $(Q, + )$, if it has the form:

$$
xy = \varphi x + c + \beta y\,\,\,(xy = \alpha x + c + \psi y),
$$
where $\beta$ (accordingly  $\alpha$) is a permutation of the set
$Q$, $\varphi \in Aut(Q, + )\,\,\,(\psi \in Aut(Q, + ))$.}

{\bf Definition 1.6.} {\it A quasigroup $(Q,\cdot)$ is called {\it
left (right) alinear} over a group $(Q, + )$, if it has the form:
$$
xy = \bar {\varphi }x + c + \beta y \quad (xy = \alpha x + c +
\bar {\psi }y),
$$
where $\beta$ (accordingly  $\alpha$) is a permutation  of the set
$Q$, $\bar {\varphi }\left( accordingly \,\,{ \bar \psi}\right)$
is an antiautomorphism of the group  $(Q, + )$.}

{\bf Definition 1.7.} {\it A quasigroup  $(Q, \cdot )$ is called
{\it of mixed type of linearity of first kind (second kind)}, if
it has the form:
$$
xy = \varphi x + c + \bar{\psi}y  \quad  (xy = \bar {\varphi}x + c
+ \psi y),
$$
where  $\varphi \in Aut(Q, + )$ ($\psi \in Aut(Q, + ))$,
$\bar{\psi}$ (accordingly $\bar {\varphi}$) is an antiautomorphism
of the group $(Q, + )$.}

All these classes shall be named classes of different type of
linearity.

An important subclass of the linear quasigroups are {\it medial
quasigroups}. A quasigroup $(Q,\cdot)$ is called {\it medial}, if
the following identity holds: $xy \cdot uv = xu \cdot yv$. By the
theorem of Bruck-Toyoda [4-6] any medial quasigroup is linear over
an abelian group with the condition $\varphi\psi=\psi\varphi$,
where $\varphi$, $\psi$ are automorphisms of the abelian group.

Medial quasigroups were  researched  by many algebraists, namely
R.H. Bruck [4], K. Toyoda [5], D.S. Murdoch [6], T. Kepka  and P.
Nemec [7,8] , K.K. Shchukin [9,10], V.A. Shcherbacov  [11] and
others, and this class plays special role in the theory of
quasigroups.

Another important subclass of linear quasigroups are  {\it
$T$-quasigroups} which were introduced and researched in detail by
T. Kepka and P. Nemec in [7,8]. According to their definition,
$T$-quasigroups are quasigroups of the form (\ref{eq1}), where
$(Q, + )$ is an abelian group and unlike medial quasigroups
$\varphi$ and $\psi$ not necessarily commute. Later  G.B.
Belyavskaya characterized the class of $T$-quasigroups by a system
of two identities [12].

\begin{center}
{\bf II.  Some  properties of linear and alinear quasigroups}
\end{center}
\quad In this part   some elementary properties of linear and
alinear quasigroups are established. Note that some other
properties of linear quasigroups are given as needed. Also we give
some necessary facts from the theory of quasigroups.


All left $(L_{a}: L_{a} x = ax)$ and right $(R_{a}: R_{a} x = xa)$
permutations of a quasigroup $(Q,\cdot)$ generate the group which
is called the { \it multiplication group} of $(Q,\cdot)$ and is
denoted by $G(\cdot)$ or $M(Q,\cdot)$.

The permutation $\alpha \in M(Q,\cdot)$ is called {\it inner} with
respect to a fixed element $h \in Q$, if $\alpha h = h$. All inner
permutations with respect to the element $h \in Q$ generate a
group which  is called the { \it group of inner permutations} of
the quasigroup $(Q,\cdot)$ and is denoted by $I_{h}$ or
$I_{h}(\cdot)$.

The action of the group $I_{h}(\cdot)$ on the quasigroup
$(Q,\cdot)$ is well known [10]. For example, if for a quasigroup
$(Q,\cdot)$ the group $I_{h}(\cdot)$ is normal subgroup in the
group $M(Q,\cdot)$, then $(Q,\cdot)$ is an abelian group  [10].
According to Theorem 4.4 from $[2]$ the group  $I_{h}(\cdot)$ is
generated by the following permutations: $R_{x,y}, L_{x,y}$ and
$T_{x}$, where $R_{x,y}=R_{x \bullet y}^{-1} R_{y} R_{x},$
$L_{x,y} = L_{x \circ y}^{-1} L_{x} L_{y},$ $T_{x}=L_{\sigma
x}^{-1} R_{x},  x \bullet y = L_{h}^{-1} (hx \cdot y), \quad x
\circ y = R_{h}^{-1} (x \cdot yh), \quad \sigma = R_{h}^{-1}
L_{h}$, $\sigma \in I_{h} (\cdot)$. If the quasigroup $(Q, \cdot)$
is a group, then $I_{h}(\cdot)=I_0(\cdot)=Int(Q, \cdot),$ i.e. the
group $I_0 (\cdot)$ is the group of inner automorphisms of the
group $(Q, \cdot)$.

Let  $(Q,\cdot)$ be a linear quasigroup:
\begin{equation}
\label{1.2.1} xy = \varphi x + c + \psi y.
\end{equation}
Here and henceforth the  $<\ldots >$  brackets will be replace the
word "generated".

{\bf Theorem  2.1.} {\it Let  $I_{h} (\cdot)$ be the group of
inner permutations of the linear quasigroup $(Q,\cdot)$. Then
$$
I_{h}(\cdot) =<R_{e_{h}}, L_{f_{h}}, T_{x}>,
$$
where $ R_{e_h} = \tilde{R}_{-\varphi h + h} \varphi, \quad
L_{f_{h}}= \tilde{L}_{h-\psi h} \psi, \quad T_{x} =
\{\tilde{L}_{h-(\psi^{-1}c +x)-\psi^{-1}\varphi h}
\tilde{R}_{\psi^{-1} c + x} \psi^{-1}\varphi,\,\,\, x \in Q\},$
and $ \tilde {R}_a x = x + a, \quad \tilde{L}_a x = a + x,$ for
all  $x \in Q$}.

\proof
According to Theorem 4.4 from  $[2],$
$$
I_{h}(\cdot)=< R_{x,y}, L_{x,y}, T_{x}>,
$$
where
$$
R_{x,y}=R_{x\bullet y}^{-1} R_{y} R_{x},\,\,\, L_{x,y}=L_{x \circ
y}^{-1} L_x L_y,\,\,\ T_{x} = L_{\sigma x}^{-1} R_{x},
$$
$$
x\bullet y =L_h^{-1} (hx \cdot y), \quad x \circ y = R_h^{-1} (x
\cdot yh), \quad \sigma = R_h^{-1} L_{h}.
$$
First we note that from  (\ref{1.2.1}) it follows that
$$
L_x = \tilde {L}_{\varphi x + c} \psi, \quad L_x^{-1} =
\psi^{-1}\tilde {L}_{\varphi x + c}^{-1}, \quad R_{y} = \tilde
{R}_{c + \psi y} \varphi, \quad R_{y}^{-1} = \varphi^{-1}\tilde
{R}_{c + \psi y}^{-1}.
$$
Then
$$
x\bullet y = L_h^{-1} (hx \cdot y)=\psi^{-1}\tilde {L}_{\varphi
h+c}^{-1} [\varphi (\varphi h + c + \psi x) + c + \psi y] = \psi
^{-1}\tilde {L}_{\varphi h + c}^{-1} [\varphi^{2}h +
$$
$$
+\varphi c + \varphi \psi x + c + \psi y] = \psi^{-1}[-c-\varphi h
+ \varphi^{2}h + \varphi c + \varphi \psi x + c + \psi y] = -
\psi^{-1}c -
$$
$$- \psi^{-1}\varphi h + \psi^{-1}\varphi^{2}h + \psi^{-1}\varphi c + \psi^{-1}\varphi
\psi x + \psi^{-1}c + y.
$$
Whence
$$
R_{x,y} (t) = R_{x\bullet y}^{-1} R_{y} R_{x} (t) = R_{x\bullet
y}^{-1} (tx \cdot y) = \varphi^{-1}\tilde {R}_{c + \psi (x\bullet
y)}^{-1} [\varphi (\varphi t + c + \psi x) + c + \psi y] =
$$
$$
= \varphi^{-1}[\varphi^{2}t + \varphi c + \varphi \psi x + c +
\psi y - \psi (x \bullet y) - c] = \varphi t + c + \psi x +
\varphi^{-1}c +
$$
$$
 + \varphi ^{ - 1}\psi y - \varphi ^{ - 1}\psi [ - \psi ^{ - 1}c - \psi ^{ -
1}\varphi h + \psi ^{ - 1}\varphi ^2h + \psi ^{ - 1}\varphi c +
\psi ^{ - 1}\varphi \psi x + \psi ^{ - 1}c +
$$
$$
+y] - \varphi^{-1}c=\varphi t + c + \psi x + \varphi^{-1}c +
\varphi^{-1}\psi y - \varphi^{-1}\psi y - \varphi^{-1}c - \psi x -
c - \varphi h + h +
$$
$$
+\varphi^{-1}c - \varphi^{-1}c = \varphi t - \varphi h + h =
\tilde {R}_{-\varphi h + h} \varphi (t),
$$
i.e. $R_{x,y} = \tilde {R}_{-\varphi h + h} \varphi$  is
independent of  $x,y$. But  $R_{e_h, e_h} = R_{e_h}$. Therefore
$$
R_{x,y} = R_{e_h}=\tilde {R}_{-\varphi h + h} \varphi.
$$
$L_{x,y} $ is  calculated similarly.

Let us now calculate  $T_{x}$:
$$
 T_{x} (t) = L_{\sigma x}^{-1} R_{x} (t) = \psi^{-1}\tilde {L}_{\varphi
(\sigma x) + c}^{-1} \tilde {R}_{c + \psi x} \varphi (t) =
\psi^{-1}(-c-\varphi(\sigma x) + \varphi t + c + \psi x)=
$$
$$ =-\psi^{-1}c-\psi^{-1}\varphi (\sigma x) + \psi^{-1}\varphi t +
 \psi^{-1}c + x = - \psi^{-1}c - \psi^{-1}\varphi R_h^{-1} L_h (x)+
$$
$$
+ \psi^{-1}\varphi t + \psi^{-1}c + x = - \psi^{-1}c -
\psi^{-1}\varphi \varphi^{-1}\tilde {R}_{c + \psi h}^{-1} \tilde
{L}_{\varphi h + c} \psi x + \psi^{-1}\varphi t + \psi^{-1}c + x =
$$
$$
 = - \psi^{-1}c - \psi^{-1}(\varphi h + c + \psi x - \psi h
- c) + \psi^{-1}\varphi t + \psi^{-1}c + x=
$$
$$
=-\psi^{-1}c + \psi^{-1}c + h- x- \psi^{-1}c - \psi^{-1}\varphi h
+ \psi^{-1}\varphi t + \psi^{-1}c + x =
$$
$$
 = h - x - \psi^{-1}c -
\psi^{-1}\varphi h + \psi^{-1}\varphi t + \psi^{-1}c + x = \tilde
{L}_{h - (\psi^{-1}c + x)- \psi^{-1}\varphi h}
 \tilde {R}_{\psi^{-1}c + x} \psi^{-1}\varphi (t),
$$
i.e.
$$
T_{x} (t) = \tilde {L}_{h - (\psi^{-1}c + x) - \psi^{-1}\varphi h}
\tilde {R}_{\psi^{-1}c + x} \psi^{-1}\varphi (t).
$$

\endproof

{\bf Theorem 2.2.} {\it Let $(Q,\cdot)$ be an  alinear quasigroup:
\begin{equation}
\label{1.2.2} xy = \bar {\varphi}x + c + \bar{\psi } y,
\end{equation}
$I_h (\cdot)$ its group of inner permutations. Then
$$
I_h (\cdot) = < R_{x,h}, L_{h,x}, T_{x}>,
$$
where
$$
R_{x,h} = \tilde {L}_{h - hx} R_{x}, \quad L_{h,x} = \tilde
{R}_{-xh + h} L_x, \quad T_{x}=\tilde {L}_{x + \bar{\psi}^{-1}c}
\tilde {R}_{\bar{\psi}^{-1}\bar{\varphi}h - \bar{\psi}^{-1}c - x +
h} \bar{\psi}^{-1}\bar{\varphi}.
$$}

\proof
From (\ref{1.2.2}) it follows:
$$
L_{x} = \tilde {L}_{\bar{\varphi}x + c} \bar{\psi}, \quad L_x^{-1}
= \bar{\psi}^{-1}\tilde {L}_{\bar{\varphi}x + c}^{-1}, \quad R_y =
\tilde {R}_{c + \bar{\psi}y} \bar{\varphi}, \quad R_y^{-1}=
\bar{\varphi}^{-1}\tilde {R}_{c + \bar{\psi}y}^{-1}.
$$
Let us calculate  $R_{x,y}$:
\begin{align*}
&x\bullet y = L_h^{-1} (hx \cdot y)= \bar{\psi}^{-1}\tilde
{L}_{\bar {\varphi }h + c}^{-1}[\bar{\varphi}(\bar{\varphi}h + c +
\bar{\psi}x) + c + \bar{\psi}y] = \\
&= \bar{\psi}^{-1}\tilde {L}_{\bar{\varphi}h + c}^{-1}
[\bar{\varphi}\bar{\psi}x + \bar{\varphi}c + \varphi^{2}h + c
+\bar{\psi}y]=\\
& =\bar{\psi}^{-1}[-c - \bar{\varphi}h + \bar{\varphi}\bar{\psi}x
 + \bar{\varphi}c + \varphi^{2}h + c + \bar{\psi}y]=\\
& = y + \bar{\psi}^{-1}c + \bar{\psi}^{-1}\varphi^{2}h + \bar{\psi
}^{-1}\bar{\varphi}c + \bar{\psi}^{-1}\bar{\varphi}\bar{\psi}x -
\bar{\psi}^{-1}\bar{\varphi}h - \bar{\psi}c.\\
\end{align*}
$$
R_{x,y} (t) = R_{x \bullet y}^{-1} R_{y} R_{x}(t) =
\varphi^{-1}\tilde {R}_{c + \bar{\psi }(x \bullet y)}^{-1}
\tilde{R}_{c + \bar{\psi}y} \bar{\varphi}\tilde {R}_{c +
\bar{\psi}x}\bar{\varphi}(t) =
$$
$$
=\bar{\varphi}^{-1}\tilde {R}_{c + \bar{\psi}(x\bullet y)}^{-1}
\tilde {R}_{c + \bar{\psi}y} \bar{\varphi}(\bar{\varphi}t + c +
\bar{\psi}x) = \bar{\varphi}^{-1}\tilde {R}_{c + \bar{\psi}(x
\bullet y)}^{-1} (\bar{\varphi}\bar{\psi}x + \bar {\varphi}c
+\varphi^{2}t + c + \bar{\psi}y)=
$$
$$
=\bar{\varphi}^{-1}(\bar{\varphi}\bar{\psi}x + \bar{\varphi}c +
\varphi^{2}t + c + \bar{\psi}y - \bar{\psi}(x \bullet y)-c) =
$$
$$
=-\bar{\varphi}^{-1}c -\bar{\varphi}^{-1}\bar{\psi}(x\bullet y) +
\bar{\varphi}^{-1}\bar{\psi}y + \bar{\varphi}^{-1}c +
\bar{\varphi}t + c + \bar{\psi}x =
$$
$$
= - \bar{\varphi}^{-1}c - \bar{\varphi}^{-1}\bar{\psi}(y +
\bar{\psi}^{-1}c + \bar{\psi}^{-1}\varphi^{2}h +
\bar{\psi}^{-1}\bar{\varphi}c+
\bar{\psi}^{-1}\bar{\varphi}\bar{\psi}x -
$$
$$
- \bar{\psi}^{-1}\bar{\varphi}h - \bar{\psi}c) +
\bar{\varphi}^{-1}\bar{\psi}y + \bar{\varphi}^{-1}c +
\bar{\varphi}t+ c + \bar{\psi}x =
$$
$$
= - \bar{\varphi}^{-1}c - \bar{\varphi}^{-1}(-c-\bar{\varphi}h +
\bar{\varphi }\bar{\psi}x + \bar{\varphi}c + \varphi^{2}h + c +
\bar {\psi }y) +
$$
$$
+ \bar {\varphi }^{ - 1}\bar {\psi }y + \bar {\varphi }^{ - 1}c +
\bar {\varphi }t + c + \bar {\psi }x = - \bar {\varphi }^{ - 1}x -
(\bar {\varphi }^{ - 1}\bar {\psi }y + \bar {\varphi }^{ - 1}c +
\bar {\varphi }h + c + \bar {\psi }x-
$$
$$
- h - \bar {\varphi}^{ -1}c) + \bar {\varphi }^{ - 1}\bar {\psi }y
+ \bar {\varphi }^{ - 1}c + \bar {\varphi }t + c + \bar {\psi}x =
- \bar {\varphi }^{ -1}c +\bar{\varphi }^{ - 1}c + h - \bar {\psi
}x - c - \bar {\varphi }h -
$$
$$
- \bar {\varphi }^{ - 1}c - \bar {\varphi }^{ - 1}\bar {\psi }y +
\bar {\varphi }^{ - 1}\bar {\psi }y + \bar {\varphi }^{ - 1}c +
\bar{\varphi }t + c + \bar {\psi }x =
$$
$$
= h - \bar {\psi }x - c - \bar {\varphi }h + \bar {\varphi }t + c
+ \bar {\psi }x = h - (\bar {\varphi }h + c + \bar {\psi }x) +
\bar {\varphi}t + c + \bar {\psi }x =
$$
$$
= h - hx + tx = \tilde {L}_{h - hx} \tilde{R}_x (t).
$$

Hence  $R_{x,y} = \tilde {L}_{h - hx} \tilde{R}_x $. Similarly is
calculated $L_{x,y}$.

Let us calculate  $T_x = L_{\sigma x}^{ - 1} R_x $, where $\sigma
= R_h^{ - 1} L_h $, considering that
\begin{align*}
&\sigma x = R_h^{ - 1} L_h x = \bar {\varphi }^{ - 1}\tilde {R}_{c
+ \bar {\psi }h}^{ - 1} \tilde {L}_{\bar {\varphi }h + c} \bar
{\psi }x = \bar {\varphi }^{ - 1}\tilde {R}_{c + \bar {\psi }h}^{
- 1} (\bar {\varphi }h + c+ \bar {\psi }x) = \\
&= \bar {\varphi }^{ - 1}(\bar {\varphi }h + c + \bar {\psi }x -
\bar {\psi }h - c) = - \bar {\varphi }^{ - 1}c - \bar {\varphi }^{
- 1}\bar {\psi }h + \bar {\varphi }^{ - 1}\bar {\psi }x + \bar
{\varphi }^{ - 1}c +h,
\end{align*}
we have
\begin{align*}
&T_x (t) = L_{\sigma x}^{ - 1} R_x (t) = \bar {\psi }^{ - 1}\tilde
{L}_{\bar {\varphi }(\sigma x) + c}^{ - 1} \tilde {R}_{c + \bar
{\psi }x} \bar {\varphi }(t) = \bar {\psi }^{ - 1}( - c - \bar
{\varphi }(\sigma x) + \bar{\varphi }t + c + \bar {\psi }x) = \\
&= x + \bar {\psi }^{ - 1}c + \bar {\psi }^{ - 1}\bar {\varphi }t
-\bar {\psi }^{ - 1}\bar {\varphi }\sigma x - \bar {\psi }^{ - 1}c = \\
&= x + \bar {\psi }^{ - 1}c + \bar {\psi }^{ - 1}\bar {\varphi }t
- \bar {\psi }^{ - 1}\bar {\varphi }( - \bar {\varphi }^{ - 1}c -
\bar {\varphi }^{ - 1}\bar {\psi }h + \bar {\varphi }^{ - 1}\bar
{\psi }x + \bar{\varphi}^{-1}c + h) - \bar{\psi}^{-1}c = \\
&= x + \bar {\psi }^{ - 1}c + \bar {\psi }^{ - 1}\bar {\varphi }t
- ( - \bar {\psi }^{ - 1}c - h + x + \bar {\psi }^{ - 1}c - \bar
{\psi}^{-1}\bar {\varphi }h) - \bar {\psi }^{ - 1}c = \\
&= x + \bar {\psi }^{ - 1}c + \bar {\psi }^{ - 1}\bar {\varphi }t
+ \bar {\psi }^{ - 1}\bar {\varphi }h - \bar {\psi }^{ - 1}c - x +
h+ \bar{\psi}^{-1}c - \bar{\psi}^{-1}c = \\
&= x + \bar {\psi }^{ - 1}c + \bar {\psi }^{ - 1}\bar {\varphi }t
+\bar {\psi }^{ - 1}\bar {\varphi }h - \bar {\psi }^{ - 1}c - x + h = \\
&=\tilde{L}_{x+\bar{\psi}^{-1}c}\tilde{R}_{\bar{\psi}^{-1}\bar{\varphi}h-\bar{\psi}^{-1}c-x+h}\bar{\psi}^{-1}\bar{\varphi}(t).
\end{align*}

Thus, $T_{x} = \tilde{L}_{x + \bar {\psi}^{-1}c} \tilde
{R}_{\bar{\psi}^{-1}\bar{\varphi}h - \bar{\psi}^{-1}c - x+
h}\bar{\psi}^{-1}\bar{\varphi}$.
\endproof

{\bf Theorem 2.3.} {\it Let  $(Q,\cdot)$ be a linear (alinear)
quasigroup of the form:
$$
xy = \varphi x + \psi y \quad (xy = \bar{\varphi }x +
\bar{\psi}y),
$$
$\varphi, \psi \in Aut(Q, + )$ ($\bar {\varphi },\bar {\psi }$ are
antiautomorphisms of $(Q,+)$).  If $(H,\cdot)$ is a subquasigroup
of the quasigroup $(Q, \cdot)$, $0 \in H$, where 0 is the zero
element of $(Q,+)$, then $(H,+)$ is a subgroup of the group
$(Q,+)$, moreover
$$
(H,\cdot )\mathop \triangleleft_ {-}(Q,\cdot ) \Leftrightarrow
(H,+)\mathop \triangleleft_{-}(Q,+).
$$}

\proof Let $(H,\cdot)$ be a subquasigroup of the quasigroup
$(Q,\cdot)$, $xy = \varphi x + \psi y$, then $x + y =
\varphi^{-1}x \cdot \psi^{-1}y=R_0^{-1} x \cdot L_0^{-1} y$. If
$x, y, 0 \in H$, then $x + y \in H$, $ - x \in H$, i.e. $(H, + )$
is a subgroup of the group $(Q,+ )$.

According to  [1], if $(H,\cdot)\mathop
\triangleleft\limits_{-}(Q,\cdot)$, then  $I_{h} H = H$ for any $h
\in Q,$ i.e. $H$ is invariant with respect to any permutation from
 $I_h $. Let  $h = 0$, then $R_{e_0 } H = R_{0} H
= \varphi H = H, \quad L_{f_0 } H = L_0 H = \psi H = H$ and
according to Theorem 2.1
$$
T_{x} = \tilde{L}_{- x} \tilde{R}_{x}\psi^{-1}\varphi.
$$
If  $(H,\cdot )\mathop \triangleleft\limits_{-}(Q,\cdot )$, then
$T_{x} H= H$ or  $- x + \psi^{-1}\varphi H  + x = H.$ Then $ - x +
H + x = H$, i.e. $(H, + )\mathop \triangleleft\limits_{-}(Q,+)$.

The converse is easily verified.

Similarly the theorem can be proved for alinear quasigroups.
\endproof

We shall denote for convenience a quasigroup as $A$ or $(Q, A)$.
It is known [1] that with  each  quasigroup $A$ the next five
quasigroups are connected:

$$
A^{-1}, \,\,\, {}^{-1}A,\,\,\, {}^{-1}(A^{-1}),\,\,
({}^{-1}A)^{-1},\,\,\,
[{}^{-1}(A^{-1})]^{-1}={}^{-1}[({}^{-1}A)^{-1}]=A^{*}
$$
These quasigroups are called {\it inverse quasigroups or
parastrophies.}

It is  known [7, 8] that T-quasigroups are invariant under
parastrophies. Below we shall establish the form of the
parastrophies of linear and alinear quasigroups.

For this we shall re-write (\ref{1.2.1}) in the form:

$$
A(x,y) = \varphi x + c + \psi y.
$$

It turns out that linear quasigroups, unlike T-quasigroups, are
not invariant under parastrophies, namely the following holds:

{\bf Proposition 2.1.} {\it Let  $(Q, A)$ be a linear quasigroup:
$$
A(x,y) = \varphi x + c + \psi y.
$$
Then
$$
A^{-1}(x,y) = \bar{\varphi}_{1} x + c_{1} + \psi_{1} y, \quad
{}^{-1}A(x,y) = \varphi_2 x + c_{2} + \bar{\psi}_2 y,
$$
$$
{}^{-1}(A^{-1})(x,y) = \varphi_{3} y + c_{3} + \bar{\psi}_{3} x,
\quad ({}^{-1}A)^{-1}(x,y) = \bar{\varphi}_{4} y + c_{4} +
\psi_{4} x,
$$
$$
A^{*}(x,y) = \varphi y + c + \psi x,
$$
where $\varphi _i ,\psi _i \in Aut(Q, + ),\ \bar {\varphi }_i
,\bar {\psi }_i $ are  antiautomorphisms of the group   $(Q, + )$,
$c_i \in Q, i=1,2,3,4.$}

\proof
Let $A(x,y)=z$. Then $A^{-1}(x,z) = y,$
$$
\begin{array}{l}
z = \varphi x + c + \psi y \Rightarrow z = \varphi x + c + \psi
A^{-1}(x,z) \Rightarrow \\
\Rightarrow \psi A^{-1}(x,z) =-c-\varphi x + z
\Rightarrow A^{-1}(x,z) = \\
= - \psi^{-1}c - \psi ^{ - 1}\varphi x + \psi^{-1}z = \\
=-\psi^{-1}c - \psi^{-1}\varphi x + \psi^{-1}c - \psi^{-1}c +
\psi^{-1}z = \bar{\varphi}_{1} x + c_{1} + \psi_{1} z, \\
 \end{array}
$$
where $\bar{\varphi}_{1} x = J\psi^{-1}c + J\psi^{-1}\varphi x +
\psi^{-1}c$ is an antiautomorphism of the group  $(Q,+)$, $c_{1} =
J\psi^{-1}c,$ $\psi_{1} = \psi^{-1}$.

Hence, $A^{-1}(x,y) = \bar{\varphi}_{1} x + c_{1} + \psi _{1} y$.

Let us denoted $A^{-1}(x,y) = B(x,y)$. Let  $B(x,y) = t$. Then
${}^{-1}B(t,y) = x$  and  $\bar{\varphi}_{1} ({}^{-1}B(t,y)) +
c_{1} + \psi_{1} y = t \Rightarrow \bar{\varphi}_{1} ({}^{ -
1}B(t,y)) = t - \psi _1 y - c_1 \Rightarrow $ $ { }^{ - 1}B(t,y) =
- \bar {\varphi }_1^{-1} c - \bar{\varphi}_1^{ - 1} \psi _1 y +
\bar {\varphi }_1^{ - 1} t =  - \bar {\varphi }_1^{ - 1} c - \bar
{\varphi }_1^{ - 1} \psi _1 y + \bar {\varphi }_{1}^ {-1} c - \bar
{\varphi }_1^{-1} c + \bar{\varphi}_1^{-1} t =  \varphi_{3} y +
c_{3} + \bar{\psi}_{3} t,$ \\ that is
$$
 {}^{-1}B(t,y) = { }^{-1}(A^{-1})(t,y) = \varphi_{3}
y + c_{3} + \bar{\psi}_{3} t,
$$
where  $\varphi_{3} y = - \bar {\varphi}_1^{-1} c - \bar{\varphi
}_1^{-1} \psi_{1} y + \bar{\varphi}_{1}^{-1} c$ is an automorphism
of the group  $(Q,+ ), \,\,  \bar{\psi}_{3} =
\bar{\varphi}_{1}^{-1} $ is an antiautomorphism of the group $(Q,
+)$,  $c_{3} = -\bar{\varphi}_1^{-1}c$.

All other parastrophies are computed similarly.
\endproof

Note that alinear quasigroups are invariant under parastrophies,
namely, we have the following:

{\bf Proposition 2.2.} {\it Let  $(Q, A)$ be an alinear
quasigroup:
$$
A(x,y) = \bar{\varphi}x + c + \bar{\psi} y.
$$
Then
$$
A^{-1}(x,y) = \bar{\psi}_{1} y + c_{1} + \bar{\varphi}_{1} x,
\quad { }^{-1}A(x,y) = \bar{\psi}_{2} y + c + \bar{\varphi}_{2} x,
$$
$$
{}^{-1}(A^{-1})(x,y) = \bar{\varphi}_{3} x + c_{3} +
\bar{\psi}_{3} y, \quad ({}^{-1}A)^{-1}(x,y) = \bar{\varphi}_{4} x
+ c_{4} + \bar{\psi}_{4} y,
$$
$$
A^{*}(x,y) = \bar{\varphi}y + c + \bar{\psi} x,
$$
where $\bar {\varphi}_{i} $, $\bar{\psi}_{i}$ are
antiautomorphisms of the group $(Q,+)$, $c_{i} \in Q, \quad i =
1,2,3,4.$}

\proof

We shall prove the condition of the Proposition for the
parastrophy $A^{-1}(x,y)$. The remaining cases are proved
similarly.

Let  $A(x,y) = z \Rightarrow A^{ - 1}(x,z) = y$. Then
$$
\begin{array}{l}
z = \bar {\varphi }x + c + \bar {\psi }A^{ - 1}(x,z) \Rightarrow
\bar {\psi }A^{ - 1}(x,z) = - c - \bar {\varphi }x + z \Rightarrow
A^{ - 1}(x,z) = \\
 = \bar {\psi}^{ - 1}z - \bar {\psi }^{-1}\bar{\varphi }x - \bar{\psi
}^{-1}c = \bar{\psi }^{-1}z - \bar {\psi}^{-1}c + \bar
{\psi}^{-1}c - \bar{\psi }^{-1}\bar {\varphi }x - \bar
{\psi}^{-1}c =\\
= \bar {\psi }_1 z + c_1 + \bar {\varphi }_1 x,
\end{array}
$$
where $\bar{\psi}_1 = \bar {\psi}^{-1}$, $\bar{\varphi}_{1} x =
\bar{\psi}^{-1}c - \bar{\psi}^{-1}\bar{\varphi}x -
\bar{\psi}^{-1}c $ are  antiautomorphisms of the group  $(Q,+)$,
$c_1 = - \bar{\psi}^{-1}c,$

that is \,\,\,  $A^{ - 1}(x,z) = \bar {\psi }_1 z + c_1 + \bar
{\varphi }_1 x.$

\endproof

Now we consider the cases when a left and a right linear (alinear)
quasigroup  are connected between each other. First let $(Q,
\cdot)$ be linear over a loop $(Q,+)$ and consider the more common
case:

{\bf Lemma 2.1.}  {\it Let  the quasigroup $(Q, \cdot)$ be left
linear over the loop $(Q, \ast)$, $xy= \varphi x \ast \beta y$ and
right linear over the loop $(Q, \circ)$, $xy = \alpha x \circ \psi
y$, where $\alpha ,\beta$ are permutations of the set $Q$,
$\varphi \in Aut(Q, \ast )$, $\psi \in Aut(Q, \circ )$. Then $(Q,
\cdot)$ is a linear quasigroup if and only if $\alpha =
{R}^{\ast}_{\beta0}\varphi, \,\, \beta={L}^{\circ}_{\alpha{e}}\psi
$, where $0$ is the identity element of the loop $(Q,\circ)$, $e$
is the identity element of the loop $(Q,\ast), {R}^{\ast}_{a}x=x
\ast a, \,\,{L}^{\circ}_{a}x=a \circ x$}.

\proof First note that a quasigroup of the form $xy = \alpha x
\circ (c \circ \psi y)$ or $xy= (\varphi x \circ c_{1}) \circ
\beta y,$ can always be converted to the form $xy = \alpha_{1} x
\circ \psi y$ or respectively $xy= \varphi x \circ \beta_{1} y$,
where $\alpha_{1} $ and $\beta_{1}$ are some permutations of  the
set $Q$. So for convenience we re-write a left or a right linear
quasigroup in the form $xy = \alpha x  \circ \psi y$ or $xy=
\varphi x \ast \beta y$ respectively. According to the condition
of Lemma 2.1, $(Q, \cdot)$ is a left and a right linear
quasigroup, that is $xy = \alpha x  \circ \psi y = \varphi x \ast
\beta y$. Putting in the last equality first $x=0$, then $y=e$, we
get $\alpha = {R}^{\ast}_{\beta0}\varphi, \,\,
\beta={L}^{\circ}_{\alpha{e}}\psi $. The reverse condition is
obvious.

\endproof

{\bf Proposition 2.3.} {\it Let  $(Q, \cdot)$ be a left and right
linear (alinear) quasigroup. Then  $(Q, \cdot)$ is a linear (an
alinear) quasigroup.}

\proof

Let the quasigroup $(Q,\cdot )$ be a left and a right linear
quasigroup. Then  $(Q, \cdot)$ has the form: $xy = \alpha x + c +
\psi y = \varphi x \oplus c_{1} \oplus \beta y,$ where $\alpha
,\beta$ are permutations of the set $Q,\varphi \in Aut(Q, \oplus
),\psi \in Aut(Q, + )$.

But then
\begin{equation}
\label{1.3.1} xy = \alpha_{1} x + \psi y = \varphi x \oplus
\beta_{1} y,
\end{equation}

where  $\alpha_{1} х = \alpha x + c, \quad \beta_{1} y = c_{1}
\oplus \beta y$, that is the groups  $(Q,+)$ and $(Q, \oplus )$
are principally  isotopic, so by the Albert's Theorem they are
isomorphic. Moreover, from the proof of this theorem it follows
that there is such an element $k \in Q$, that  $R_{k} (x \oplus y)
= R_{k} x + R_{k} y$, where $R_{k} x = x + k$. Considering this in
(\ref{1.3.1}), we find that
$$
R_{k} (\alpha_{1} x + \psi y) =R_{k}\varphi x \oplus R_{k}\beta y,
R_{k} (\alpha_{1} x + \psi y) =R_{k} \varphi x + R_{k} \beta_{1}
y,
$$
or
$$
\alpha_1 R_{k}^{-1} x + \psi y = R_{k} \varphi R_{k}^{-1} x +
\beta_{1} y.
$$

Put in the last equality  $y =0:$ $\alpha_1 R_{k}^{-1} x =
\varphi_{1} x + d$, where $d = \beta_{1} 0$ is some element of the
set  $Q$,  $\varphi_{1} = R_{k} \varphi R_{k}^{-1} $ is
automorphism of the group $(Q,+)$.

Indeed, since  $R_{k}(x\oplus y)=R_{k}x+R_{k}y,$ then
\begin{align*}
&R_{k}\varphi R_{k}^{-1}(x+y)=R_{k}\varphi
R_{k}^{-1}[R_{k}(R_{k}^{-1}x\oplus R_{k}^{-1}y)]=\\
&=R_{k}\varphi(R_{k}^{-1}x\oplus R_{k}^{-1}y)=R_{k}[\varphi
R_{k}^{-1}x \oplus \varphi R_{k}^{-1}y]=\\
&=R_{k}[R_{k}^{-1}(R_{k}\varphi R_{k}^{-1}x+R_{k}\varphi
R_{k}^{-1}y)]=R_{k}\varphi R_{k}^{-1}x+R_{k}\varphi R_{k}^{-1}y.
\end{align*}

Hence  $R_{k}\varphi R_{k}^{-1}\in Aut(Q,+)$. According to Lemma
2.5 from [2] $\alpha _1 R_{k}^{-1}$, and thus $\alpha_{1}$ are
quasiautomorphisms of this group. Therefore $\alpha_{1} x =
\varphi_{2} x + s, \quad s \in Q,$  and   $xy = \varphi_{2} x + c
+ \psi y$, where $\varphi_2, \psi \in Aut(Q, + )$, that is $(Q,
\cdot)$ is a linear quasigroup.

Similarly, Proposition 2.3 can be proved for alinear quasigroups.
\endproof

{\bf Corollary 2.1}. {\it A left (right) linear  quasigroup
$(Q,\cdot )$ $xy = \varphi x + c + \beta y$ $(xy = \alpha x + c +
\psi y)$ is a right (left) linear quasigroup if and only if the
permutation $\beta$ (accordingly $\alpha$) is a quasiautomorphism
of the group $(Q,+)$.}

{\bf Corollary 2.2}. {\it A left (right) alinear  quasigroup
$(Q,\cdot )$ $xy = \bar{\varphi}x + c + \beta y\,\,(xy = \alpha x+
c + \bar {\psi }y)$ is a right (left) alinear quasigroup if and
only if the permutation $\beta$ (accordingly $\alpha$)  is an
antiquasiautomorphism of the group $(Q,+)$.}

{\bf Proposition 2.4.} {\it A quasigroup which is left (right)
linear and left (right) alinear is a T-quasigroup.}


\proof

Let $(Q,\cdot )$ be a left linear and a left alinear quasigroup:
$$
xy = \varphi x + c + \beta y = \bar {\varphi }x \oplus s \oplus
\gamma y,
$$
where  $\varphi \in Aut(Q, +), \bar{\varphi}$ is an
antiautomorphism of the group  $(Q, \oplus )$, $\beta$, $\gamma $
are permutations of the set $Q$. Then
\begin{equation}
\label{1.3.2} xy = \varphi x + \beta_{1} y = \bar{\varphi}x \oplus
\gamma_{1}y,
\end{equation}
where $\beta y = c + \beta y, \quad \gamma_{1} y = s \oplus \gamma
y$. Then according to Albert's Theorem $\,\,(Q, + ) \cong (Q,
\oplus )$, besides that there is an element  $k \in Q$ such that
\begin{equation}
\label{1.3.3} R_{k} (x \oplus y) = R_{k} x + R_{k} y,
\end{equation}
where  $R_k x = x + k$. Considering this in (\ref{1.3.2}), we
obtain
$$
R_{k} (\varphi x + \beta_{1} y) = R_{k} \bar{\varphi}x + R_{k}
\gamma_1 y
$$
or
$$
\varphi x + \beta_{1} y = \bar{\varphi}x + k + \gamma_{1} y.
$$
Put  $y = 0: \quad \varphi x + \beta_{1} 0 = \bar{\varphi}x + k +
\gamma_{1} 0$. Then  $\varphi x + \beta_{1} 0 - \gamma_{1} 0 - k =
\bar {\varphi}x$, $\bar{\varphi}x = \varphi x + p = R_{p} \varphi
x$, where $p = \beta _1 0 - \gamma_{1} 0 - k $ is some element of
the set $Q$.

Considering  (\ref{1.3.3}), note that  $R_{k} \bar {\varphi
}R_{k}^{-1}$ is an antiautomorphism of the group  $(Q,+)$. But
$R_{k} \bar{\varphi }R_{k}^{-1} = R_{k} R_{p} \varphi R_{k}^{-1}$.
Hence,
$$
\begin{array}{l}
R_{k} R_{p} \varphi R_{k}^{-1} (x +y) = R_{k} R_{p} \varphi
R_{k}^{-1} y + R_{k} R_{p} \varphi R_{k}^{-1} x,\\
\varphi (x + y - k) + p + k = \varphi (y - k) + p + k + \varphi(x
-k) + p + k,\\
\varphi (x + y - k) = \varphi (y - k) + p + k + \varphi (x - k),\\
\varphi x+ \varphi y - \varphi k = \varphi y - \varphi k + p + k +
\varphi x - \varphi k,\\
\varphi x + \varphi y = \varphi y - \varphi k + p + k + \varphi
x.\\
\end{array}
$$
Put $x = y = 0$: $0 = - \varphi k + p + k$. Then $\varphi x +
\varphi y = \varphi y + \varphi x$ or  $x + y = y + x$, that is
the group  $(Q,+)$ is abelian. Hence, $(Q,\cdot)$ is a left
$T$-quasigroup.

The right linear and right alinear case is proved similarly.
\endproof

Finally observe that if  $(Q,\cdot)$  is a linear and an alinear
quasigroup at the same time, then it is a $T$-quasigroup and from
Proposition 2.3 it follows that a left and right $T$-quasigroup
 is also a $T$-quasigroup.

\begin{center}
  {\bf III. Endotopies and endomorphisms of linear quasigroups}
\end{center}

{\bf Theorem 3.1.} {\it The semigroups of endomorphisms of
parastrophic quasigroups coincide:
\begin{equation}\label{1.6.1}
End\left( {Q,\cdot}\right) = End\left( {Q,
{}^\sigma(\cdot)}\right),
\end{equation}
where  $\left(Q,\cdot\right)$ is a quasigroup, $(Q,
{}^{\sigma}\left(\cdot\right))$ is its parastrophy, $End\left(
{Q,\cdot} \right)$ is the semigroup of  endomorphisms of the
quasigroup $(Q,\cdot)$ and $End\left(Q, {}^\sigma(\cdot)\right)$
is the semigroup of endomorphisms of its parastrophy $\left(Q,
{}^\sigma(\cdot)\right)$.}

\proof

Let $\gamma \in End\left( {Q, \cdot } \right):$ $\gamma(xy)=\gamma
x\cdot \gamma y.$ Denote $xy=z$. Then $x=z/y$. Hence $\gamma
z=\gamma(z/y)\cdot \gamma y$ or $\gamma(z/y)=\gamma z/\gamma y$.

Consequently, $\gamma \in End\left( {Q,{}^\sigma\left( \cdot
\right)} \right)$, so  $End\left( {Q, \cdot } \right) \subseteq
End\left( {Q,{}^\sigma\left( \cdot \right)} \right).$

The converse statement is obvious. Then $End\left( {Q, \cdot }
\right) = End\left( {Q,^{\sigma}\left( \cdot \right)} \right)$.

For other parastrophies the equality (\ref{1.6.1}) is proved
similarly.
\endproof

{\bf Corollary 3.1.} {\it The groups of automorphisms of
parastrophic quasigroups  coincide:}
$$
Aut(Q,\cdot)=Aut(Q, {}^\sigma(\cdot)).
$$

Denote by  $S (Q)$ the semigroup of all the transformations of the
set $Q$. Let $\left( {Q, \cdot } \right)$ be a quasigroup and
$\varphi \in End\left( {Q, \cdot } \right)$. Then
$$
\varphi(xy)=\varphi x\cdot \varphi y,
$$
\begin{equation}
\label{1.6.2} \varphi L_x = L_{\varphi x} \varphi ,
\end{equation}
\begin{equation}
\label{1.6.3} \varphi R_x = R_{\varphi x} \varphi.
\end{equation}

The following theorem is a generalization of Theorem 1 from [13].

{\bf Theorem 3.2.} {\it Let  $\left( {Q, \cdot } \right)$ be a
quasigroup and $\varphi \in S\left( Q \right)$ is such that
$\varphi L_x = L_z \varphi $, for any  $x \in Q$ and some  $z$,
depending on $x$. Then $z = \varphi x$ (that is $\varphi \in
End\left( {Q, \cdot } \right)$) if and only if $\varphi e =
e\varphi $, where $e$ is a map such that $e: a \to e(a)$ and
$a\cdot e(a)=a.$}

The Theorem is proved by analogy with the proof of the  Theorem 1
from [13].

{\bf Corollary 3.2.} {\it Let  $i$ be the right identity  element
of a quasigroup $\left(Q, \cdot \right)$ and $\varphi$ be a map of
the set $Q$ into itself such that $\varphi L_x = L_z \varphi $,
for any $x$ and some $z$, depending on $x \in Q$. Then $z =
\varphi x$ (that is  $\varphi \in End\left(Q, \cdot \right)$) if
and only if $\varphi \left( i \right) = i$.}

Now let  $\varphi \in S\left( Q \right)$ and  $\varphi \left( {x
y} \right) = \varphi y \cdot \varphi x $, that is  $\varphi$ is an
antiendomorphism of a quasigroup $\left( {Q, \cdot } \right)$.
Then in this case instead of (\ref{1.6.2}) and (\ref{1.6.3}) we
must consider the following relations:
\begin{equation}
\label{1.6.4} \varphi L_x = R_{\varphi x} \varphi ,
\end{equation}
\begin{equation}
\label{1.6.5} \varphi R_x = L_{\varphi x} \varphi.
\end{equation}

{\bf Theorem 3.3.} {\it Let  $\left( {Q, \cdot } \right)$ be a
quasigroup and $\varphi \in S\left( Q \right)$ is such that
$\varphi L_x = R_z \varphi $, for any  $x$ and some  $z$,
depending on $x \in Q$. Then $z = \varphi x $ (that is $\varphi$
is an antiendomorphism  of $\left( {Q, \cdot } \right)$) if and
only if $\varphi e = f\varphi $, where $e$ is a map such that $e:
a \to e(a)$ and $a\cdot e(a)=a$ and $f$ is a map such that $f:$ $a
\to f\left( a \right)$ and $f\left( a \right) \cdot a = a$.}

\proof

Since $\varphi L_x = R_z \varphi$, we have that $\varphi L_x y =
R_z \varphi y$, i.e. $\varphi \left( {xy} \right) = \varphi y
\cdot z$. In particular, for  $y = e\left( x \right)$: $\varphi
\left( {x \cdot e\left( x \right)} \right) = \varphi e\left( x
\right) \cdot z$, whence $\varphi \left( x \right) = \varphi
e\left( x \right) \cdot z$.

But $\varphi (x) = f\varphi (x) \cdot \varphi (x)$. Then $\varphi
e(x)\cdot z = f\varphi (x) \cdot \varphi (x)$.

Assume now that $\varphi e = f\varphi $. After reduction we get $z
= \varphi (x)$, that is  $\varphi$ is an antiendomorphism of the
quasigroup $(Q, \cdot )$.

Conversely, for an antiendomorphism $\varphi$ the following holds:
$$
f \varphi (x) \cdot \varphi (x) = \varphi (x) = \varphi (x \cdot
e(x)) = \varphi e(x) \cdot \varphi (x)
$$
hence  $f \varphi (x) = \varphi e(x)$ for all $x \in Q$.
\endproof

{\bf Corollary 3.3.} {\it Let  $(Q, \cdot)$ be a loop with the
identity element  $j$ and  $\varphi \in S(Q)$ be such that
$\varphi L_{x} = R_{z} \varphi $ for any  $x$ and some $z,$
depending on $x \in Q$. Then $z = \varphi(x)$ (that is $\varphi$
is an antiendomorphism of $(Q, \cdot)$)  if and only if
$\varphi(j)=j$.}

\proof

In a loop  $(Q,\cdot)$ we have $e(x)=j = f(x)$  for any $x \in Q$.
Let $\varphi$ be an  antiendomorphism of  the loop  $(Q, \cdot)$.
Then
$$
\varphi(j) = \varphi (e(x)) = \varphi (f(x))=f\varphi(x) = j.
$$
Conversely,  if  $\varphi L_{x} = R_{z} \varphi $ and $\varphi (j)
= j$, then since  $j = e(x) = f(x)$, we have $\varphi (e(x)) =
e(x) = j = f(x) = f\varphi (x)$, whence $\varphi e = f\varphi $.
Then according to Theorem 3.3  $\varphi$ is an antiendomorphism of
the loop $(Q, \cdot)$.
\endproof

Recall that the triplet $T = (\alpha ,\beta ,\gamma )$ of maps of
a quasigroup $(Q,\cdot)$ into itself is called an \emph{endotopy}
of the quasigroup $(Q,\cdot)$, if the identity:
$$
\gamma (xy) = \alpha x \cdot \beta y,
$$
is true for any $x,y \in Q$.

In case when  $\alpha = \beta = \gamma $, the triplet $T =
(\gamma,\gamma, \gamma)$ is called an \emph{endomorphism} of the
quasigroup $(Q,\cdot )$.

Obviously the set of all endotopies of the quasigroup $(Q,\cdot)$
is a semigroup with identity. Let us denote this semigroup by
$Ent(Q, \cdot)$.

{\bf Theorem 3.4.} {\it If the quasigroups  $(Q,\cdot )$ and
$(Q,\circ)$ are isotopic: $\gamma (x \circ y) = \alpha x \cdot
\beta y,\,\,\, (\circ)=(\cdot) T,\,\,\, T=(\alpha ,\beta
,\gamma),$ then their semigroups of endotopies are conjugate:}
\begin{equation}\label{1.6.6}
Ent(Q, \cdot )=T^{-1}Ent(Q,\circ)T.
\end{equation}

\proof

Indeed, let $(\circ)=(\cdot)T$ and $S \in Ent(Q, \circ )$, that is
$( \circ )S = ( \circ )$. Then $(( \cdot )T)S = ( \cdot )T$,
therefore $( \cdot )TS = ( \cdot )T$ and  $( \cdot )TST^{ - 1} = (
\cdot )$,
so $T Ent(Q, \circ )T^{ - 1} \subseteq Ent(Q, \cdot )$, whence \\
$Ent(Q, \circ ) \subseteq T^{-1}Ent(Q, \cdot )T$.

On the other side, $( \circ )T^{-1} = (\cdot)$. Let $S^{'} \in
Ent(Q, \cdot )$, that is  $(\cdot )S^{'} = ( \cdot )$. Then $(
\circ )T^{ - 1}S^{'} = ( \circ )T^{ - 1}$ and $( \circ
)T^{-1}S^{'} = ( \circ )$, whence  $T^{-1}Ent(Q, \cdot )T
\subseteq Ent(Q, \circ )$. Considering the reverse inclusion we
get  (\ref{1.6.6}).
\endproof

{\bf Corollary 3.4. [1]}. {\it If the quasigroups $(Q, \cdot)$ and
$(Q, \circ)$ are isotopic then their groups of automorphisms are
isomorphic, namely}
$$
A \nu t(Q, \cdot) = T^{-1}A \nu t(Q,\circ)T.
$$

Note that the third component $\gamma$ of an endotopy $T=(\alpha,
\beta, \gamma)$ is called \emph{a quasiendomorphism}  of the
quasigroup $(Q, \cdot )$. The structure of the quasiautomorphisms
of a group is well known (see [1]). As in the case of
quasiautomorphisms, the quasiendomorphisms of groups have a simple
structure, namely:

{\bf Proposition 3.1.} {\it Any quasiendomorphism  $\gamma $ of a
group $(Q, + )$ has the form:
\begin{equation}
\label{1.6.7} \gamma = \tilde {R}_s \gamma_ \circ,
\end{equation}
where $\gamma _ \circ \in End(Q, + )$, $s \in Q$, and reversely,
the map $\gamma $, defined by the equality (\ref{1.6.7}), is a
quasiendomorphism of the group $(Q, + )$.}

\proof

Let $\gamma$ be a quasiendomorphism of a group $(Q,+)$, that is
there is an endotopy  $T = (\alpha,\beta,\gamma): \quad \gamma (x
+ y) = \alpha x + \beta y$. Put  $x = 0$ in the last equality,
where $0$ is the identity element of the group $(Q, + )$. Then
$\gamma y = \alpha 0 + \beta y = l + \beta y = \tilde {L}_l \beta
y$, $\gamma = \tilde {L}_l \beta $, $\beta = \tilde {L}_l^{-1}
\gamma $, $l = \alpha 0$. Similarly taking $y = 0$ we have:
$\gamma х = \alpha x + \beta 0 = \alpha x + k = \tilde {R}_{k}
\alpha х$, $\gamma = \tilde {R}_{k} \alpha $, $\alpha = \tilde
{R}_{k}^{- 1} \gamma $, $k = \beta 0$. But in  any group $\tilde
{L}_l^{-1} = \tilde {L}_{-l} $, $\tilde {R}_{k}^{-1} =
\tilde{R}_{-k}$.

Hence,
$$
\gamma (x + y) = \alpha x + \beta y = \tilde {R}_{-k} \gamma x +
\tilde {L}_{-l} \gamma y = (\gamma x + (-k)) + ((-l) + \gamma y),
$$
\begin{equation}
\label{1.6.8}\gamma (x + y) = \gamma x + ( - \gamma 0) + \gamma y,
\end{equation}
since $(-k) + (-l) = - (l + k)=-(\alpha 0 + \beta 0) = - \gamma
0$. Let  $\gamma x + ( - \gamma 0) = \gamma _0 x$. Adding
$(-\gamma 0)$ to the both parts of (\ref{1.6.8}), we have $\gamma
(x + y) + ( - \gamma 0) = \gamma x + ( - \gamma 0) + \gamma y + (
- \gamma 0)$, $\gamma _0 (x + y) = \gamma _0 x + \gamma _0 y$,
that is $\gamma _0 \in End(Q, + )$. From $\gamma x + (-\gamma 0) =
\gamma_0 x$ follows (\ref{1.6.7}).

Conversely, any endomorphism and any transformation of a group are
quasiendomorphisms, so the product  $\tilde {R}_s \gamma _0$ is
also a quasiendomorphism of that group.
\endproof

{\bf Corollary 3.5.} {\it Let  $\gamma $  be a quasiendomorphism
of a group $(Q, + )$. Then
$$
\gamma \in End(Q, + ) \Leftrightarrow \gamma 0 = 0,
$$
where $0$ is the identity element of the group  $(Q,+)$.}

Indeed, let $\gamma $ be a quasiendomorphism and  $\gamma 0 = 0$.
Then $\gamma 0 = \tilde {R}_s \gamma_0 0 = R_s 0 = s$ and $s = 0$.
Hence, $\gamma = \gamma_0 $.

As noted by I.A. Golovko in [14], any endotopy of a group  $(Q, +
)$ has the form:
\begin{equation}
\label{1.6.9} T = (\tilde {L}_{a},\,\, \tilde {R}_{b}, \,\,\tilde
{L}_{a} \tilde {R}_{b} )\theta ,
\end{equation}
where $\theta \in End(Q, + )$, $a, b \in Q$.

Using  (\ref{1.6.6}) and  (\ref{1.6.9}), it is easy to find a
common form of any endotopy of a linear (alinear) quasigroup $(Q,
\cdot )$.

Indeed, let  $(Q,\cdot)$ be a linear quasigroup: $xy = \varphi x +
c + \psi y = \tilde {R}_c \varphi x + \psi y$. Then
$(\cdot)=(+)(\tilde {R}_c \varphi ,\psi ,\varepsilon )$. If  $S
\in Ent(Q, + )$, then by  (\ref{1.6.9}) $S =$ $(\tilde
{L}_{a},\,\,\tilde {R}_{b}, \,\,\tilde {L}_{a} \tilde
{R}_{b})\theta $, where $\theta \in Ent(Q, + )$, $a, b \in Q$.
Then according to  (\ref{1.6.6}) for any  $P \in Ent(Q, \cdot )$
we have:
$$
\begin{array}{l}
P = TST^{ - 1} = (\tilde {R}_c \varphi, \psi, \varepsilon )
(\tilde {L}_{a} \theta, \tilde {R}_{b} \theta, \tilde {L}_{a}
\tilde {R}_{b} \theta)(\varphi^{-1}R_c^{-1},
\psi^{-1},\varepsilon) =\\
=(\tilde {R}_c \varphi \tilde {L}_{a} \theta \varphi^{-1}\tilde
{R}_c^{-1}, \psi \tilde {R}_{b} \theta \psi^{-1}, \tilde {L}_{a}
R_{b} \theta) = (\tilde {R}_c \tilde {L}_{\varphi a} \varphi
\theta \varphi^{-1}\tilde {R}_c,\,\, R_{\psi b} \psi \theta
\psi^{-1},\tilde {L}_{a} \tilde {R}_{b} \theta),
\end{array}
$$
where $\varphi, \psi \in Aut(Q, + ), \theta \in End(Q, + ),\,\, a,
b, c \in Q.$

So, have proved:

{\bf Theorem 3.5.} {\it  Any endotopy of a linear quasigroup  $(Q,
\cdot )$:  $xy = \varphi x + c + \psi y$ has the form:}
\begin{equation}
\label{1.6.10} {\rm P} = (\tilde {R}_c \tilde {L}_{\varphi a }
\varphi \theta \varphi^{-1}\tilde {R}_{- c}, \tilde {R}_{\psi b}
\psi \theta \psi^{-1}, \tilde {L}_{a} \tilde {R}_{b} \theta ),
\end{equation}
{\it where} $\varphi,\psi \in Aut(Q, + )$, $\theta \in End(Q, +
)$, $a, b, c \in Q$.

Similarly can be showed that any endotopy an alinear
quasigroup$(Q, \cdot )$: $xy = \tilde {\varphi}x + c + \tilde
{\psi}y$ has the form:
\begin{equation}
\label{1.6.11} \overline P = \left( {\tilde {R}_d \overline
\varphi \theta \overline \varphi^{-1}\tilde {R}_{-c}, \tilde
{L}_{\overline \psi b} \overline \psi \theta \overline \psi ^{ -
1},\tilde {L}_{a} \tilde {R}_{b} \theta } \right).
\end{equation}

{\bf Corollary 3.6.} {\it Any endotopy of a T-quasigroup (medial
quasigroup)  $(Q,\cdot )$: $xy= \varphi x+c+\psi y$ has the form:
\begin{equation}
\label{1.6.12} P=(R_{d} \varphi \theta \varphi{}^{-1}R_{-c},\,\,
\overline{R}_{\psi b} \psi\theta \psi^{-1}, \,\,
\overline{L}_{a+b}\theta),
\end{equation}
where $d=c+\varphi a, \,\, d \in Q, \theta \in End(Q , +).$}

{\bf Corollary 3.7.} {\it Any endotopy a quasigroup of mixed type
of linearity of first kind:   $xy=\varphi x+c+\bar{\psi}y$
(accordingly linearity of second kind:  $xy =\bar{\varphi}x+c+\psi
y$) has the form:
$$
P=\left(\tilde{R}_{c}\tilde{L}_{\varphi a} \varphi
\theta\varphi{}^{-1}\tilde{R}_{-c}, \,\,
\bar{\psi}\tilde{R}_{b}\theta \bar{\psi}^{-1},\,\,
\tilde{L}_{a}\tilde{R}_{b}\theta\right),
$$
$$
\left(P=\left(\tilde{R}_{c}\tilde{L}_{\varphi
a}\bar{\varphi}\theta\bar{\varphi}{}^{-1}\tilde{R}_{-c},\,\,\,
\psi \tilde{R}_{b}\theta \bar{\psi}^{-1},\,\,\,
\tilde{L}_{a}\tilde{R}_{b}\theta\right)\right),
$$
where $\varphi, \psi \in Aut(Q,+),$ $\bar{\varphi}, \bar{\psi}$
are antiautomorphisms of the group $(Q,+),$ $a, b, c$ are fixed
elements of $Q,$ $\theta \in End (Q, +),$ $\tilde{R}_{a},
\tilde{L}_{b}$ are transformations of the group $(Q, +):$
$\tilde{R}_{a}x =x+a,$ $\tilde{L}_{a}x=a+x.$}

{\bf Corollary 3.8.} {\it Any endomorphism  $\gamma$ of a linear
quasigroup $\left( {Q, \cdot}\right): xy = \varphi x + c + \psi y$
can be presented in the form:
$$
\gamma = \tilde{R}_{c}\tilde{L}_{\varphi a}\varphi \theta
\varphi{}^{-1}\tilde{R}_{-c}=\tilde{R}_{\psi b}\psi \theta
\psi^{-1}=\tilde{L}_{a}\tilde{R}_{b}\theta,
$$

 where $\theta\in End(Q,+),$ $a \in Q.$}

\proof It is follows from the coincidence of all components of the
endotopy of a linear quasigroup.
\endproof

{\bf Remark.} Similarly, any endomorphism (automorphism) of
alinear quasigroups,  quasigroups of mixed type of linearity,
T-quasigroups and medial quasigroups can be presented by
automorphisms, endomorphisms and transformations of respective
groups [15].
\\

{\bf Acknowledgments}

I am indebted to professors G.B. Belyavskaya  and V.A. Shcherbacov
for introducing me to this subject. Professor G. Morosanu  gave me
useful advice on writing this article. I express my gratitude to
all them.


\begin{center}
{\bf References:}
\end{center}

\begin{enumerate}
\item V.D. Belousov. {\it Foundations of the theory of quasigroups and
loops.} Moscow, Nauka, 1967.(in Russian).
\item V.D. Belousov.  {\it Balanced identities on quasigroups.} Mat. Sb.,
vol.70(112): no.1, 1966, p.55-97,(in Russian).
\item G.B. Belyavskaya,  A.Kh. Tabarov.  {\it Characterization of linear and
alinear quasigroups.} Diskretnaya matematika, v.4, no.2, 1992,
p.142-147.(in Russian).
\item R.H. Bruck.  {\it Some results in theory of
quasigroups.} - Trans. Amer. Math. Soc., 1944, vol. 55, p.19-52.
\item K. Toyoda.  {\it On axsioms of linear functions.} - Proc. Imp.
Acad.Tokyo., 1941, vol. 17, p.221-227.
\item D.S. Murdoch.  {\it Quasigroups
which satisfy certain generalized associative laws.} - Amer.
J.Math., 1939, vol. 61, p.509-522.
\item T. Kepka   and P. Nemec.  {\it $T$-quasigroups. I.} - Acta univ.
Carolin. Math.Phys., 1971, vol.12, No. 1, р.31-39.
\item T. Kepka   and P. Nemec.  {\it $T$-quasigroups.II.} - Acta univ.
Carolin. Math.Phys., 1971, vol. 12, No. 2, р.39-49.
\item K.K. Shchukin.  {\it On simple medial quasigroups.} - Math. issled. Kishinev,
1991, vol. 120, p.114-117,(in Russian).
\item K.K. Shchukin.   {\it Action a
group on a quasigroup.} Kishinev State University Printing House,
Kishinev, 1985, (in Russian).
\item V.A. Shcherbacov.  {\it On linear and inverse quasigroups and their
applications in code theory.} Dissertation of Doctor of Science,
247 pages, Chisinau, 2008.
\item  G.B. Belyavskaya. {\it Abelian
quasigroups are T-quasigroups.} Quasigroups and Related Systems,
1994, vol.1, no.1, p.1-7.
\item L.V. Safonova, K.K. Shchukin.  {\it Calculation of automorphisms and
antiautomorphisms of quasigroups.} Izvest. AS SSR Moldova,
Mathematica, 1990, vol. 3, No. 3, p.49-55,(in Russian).
\item  I.A. Golovko.  {\it Endotopies in quasigroups.} Abstracts. First
Allunion symposium on theory quasigroups and its applications.
Sukhumi,1968, p.14-15, (in Russian).
\item A.Kh. Tabarov.  {\it Homomorphisms and endomorphisms of linear and
alinear quasigroups.} Discrete Mathematics and Applications, 17
(2007), no.3, p. 253-260.

\end{enumerate}

\end{document}